\documentclass[12pt,reqno]{amsart}
\usepackage{amsfonts,amsmath,amssymb,amsthm}
\usepackage{enumitem}
\usepackage{amsaddr}
\usepackage{mathrsfs}
\usepackage{mathabx}
\usepackage{bbm}
\usepackage{cite}
\usepackage{dsfont}
\newtheorem{Theorem}{Theorem}[section]

\newtheorem{Corollary}{Corollary}[section]
\newtheorem{Remark}{Remark}[section]
\newtheorem{Definition}{Definition}[section]
\numberwithin{equation}{section}

\def \vu{\textbf{u}}

\newcommand\vN{{\bf  \nabla}}

\usepackage[a4paper,margin=1.8cm]{geometry}
\begin{document}

\sloppy
\title[]
{What means ``slightly supercritical"\\ in
Navier-Stokes equations?}
\author[]{Pedro Gabriel Fernández Dalgo}
\address{ \tiny{Escuela de Ciencias Físicas y Matemáticas, Universidad de Las Américas, Vía a Nayón, C.P.170124, Quito, Ecuador}}
\email{pedro.fernandez.dalgo@udla.edu.ec}

\date{September, 2021} 

\begin{abstract}
The goal of this paper is to clarify some differences between the critical Lebesgue space $L^3$ and the critical weak Lebesgue space $L^{3, \infty}$, when these are considered in the hypothesis of classical statements for the 3D homogeneous incompressible Navier-Stokes equations.
\end{abstract}


\maketitle
\vspace{-1cm}
\section{Introduction}

Recent interesting technical progress in the analysis of the velocities solving the Navier-Stokes equations appear associated to critical hypothesis. Let us fix some notions. The homogeneous incompressible Navier-Stokes equations
\begin{equation}\label{NS}
(\text{NS}) \left\{ \begin{array}{ll}\vspace{2mm} 
\partial_t \vu = \Delta \vu - (\vu \cdot \nabla) \vu - \nabla p \\
\vN \cdot  \vu=0, \phantom{space}
\vu(0,\cdot)=\vu_0, 
\end{array}
\right.     
\end{equation}
satisfy the following scaling property : If $({\bf u},p)$ is a solution of (NS) on $\mathbb R ^3 \times [0,T)$ with initial velocity $\vu_0$, then defining
\begin{equation*}
    {\bf u}_\lambda(x,t):= \lambda {\bf u}(\lambda x, \lambda^2 t) ,\phantom{space} p_\lambda(x,t):= \lambda^2 p(\lambda x, \lambda^2 t)  \phantom{spa} \text{and} \phantom{spa} \vu_{0, \lambda}=\lambda \vu_0 (\lambda \, \cdot),
\end{equation*}
we obtain $({\bf u}_\lambda,p_\lambda)$ is a solution of (NS) on $\mathbb R ^3 \times [0,T/\lambda^2)$ with initial data $\vu_{0, \lambda}$.

A norm or quasi-norm $ \|\cdot\|_{X_T} $, on the time-dependent vector fields 
with domain $\mathbb R ^3 \times [0,T)$, is called critical (for the velocity) if for all time-dependent vector field $\vu$ on $\mathbb R ^3 \times [0,T)$, for all $\lambda>1$,
$\|{\bf u}_\lambda \, \mathds 1 _{\mathbb R ^3 \times [0,T/\lambda^2)} \|_{X_{T}} = \|{\bf u}\|_{X_T}.$
Restricting our attention to the initial data, we say that a norm or semi-norm $\|\cdot\|_{X_0}$ defined on the vector fields with domain $\mathbb R ^3$, is critical if for all vector field ${\bf u}_{0}$ on $\mathbb R ^3$, for all $\lambda>0$,
$\|{\bf u}_{0,\lambda}\|_{X_0} = \|{\bf u}_0\|_{X_0}.$

\noindent
The most classical results concerning strong solutions of (NS) are limited by critical hypothesis (we refer to \cite{FK_64}, \cite{S62} and \cite{La67}). For example, the local existence of mild solutions with initial data in the critical Lebesgue space $L^3$ is known (see \cite{FK_64}, \cite{LR16} and \cite{LRibe}), but in the case of initial data belonging to the critical weak Lebesgue space $L^{3, \infty}$, smallness conditions are required:
\begin{equation*}
    \|\vu\|_{L^{3, \infty}(\mathbb R ^3)}:= \sup_{\alpha>0} \alpha \left[ \, \mu \left( \left\{ \, x \in \mathbb R ^3 \, : \, |\vu(x)|>\alpha \,     \right\} \right) \, \right]^{1/3}< \epsilon_{\rm univ}, \phantom{spa}\text{being $\mu$ the Lebesgue measure}.
\end{equation*}
Regularity of Leray solutions belonging to $L^{\infty}((0,T)L^{3})$ has been proved (see \cite{ESS}, \cite{Ser12}, \cite{Serbook14}, \cite{Tsaibook17}), however under the ``slightly" weaker hypothesis $L^{\infty}((0,T)L^{3,\infty})$, regularity is an open problem when the solution is not axisymmetric (we refer to\cite{BFP}).

In attempting to deal with supercritical hypotheses, the use of quantitative arguments to study concentration of some frequencies in presence of a singularity,  permits to precise a blow-up rate slightly supercritical ``with respect to the $L^3$ norm". The idea is due to Tao in \cite{Tao19}, who proved that for a finite energy solution ${\bf u}$ with viscosity $\nu=1$, which first blows up at $T^*$, we have for some absolute constant $c > 0$,
\begin{equation}
\label{tao}
    \limsup_{t \uparrow T^* } \frac {\| {\bf u} (\cdot ,t) \|_{L^3(\mathbb{R}^3)}} { \left( \text{log log log}  \frac{1} {T^* -t} \right)^c }  = + \infty.
\end{equation}
A quantitative strategy has been established in the local setting, considering type I conditions $\| {\bf u} \|_{L^\infty ((0,T^*), L^{3, \infty  })} < + \infty $,  by Barker and Prange in \cite{TB_CP21}, and also by Kang, Miura and Tsai in \cite{KMT_21}. Barker and Prange demonstrate under type I conditions that
\begin{equation*}
    \| {\bf u} (\cdot ,t) \|_{L^3(B_{x_0}(\delta))} \geq C_M \text{log} \left( \frac{1} {T^*-t} \right), \phantom{space} \text{where } \delta=  O((T^*-t)^{\frac{1}{2}-}),
\end{equation*} 
based on localized smoothing estimates in \cite{TB_CP20} while Kang, Miura and Tsai propose a direct approach in paper \cite{KMT_21}. 
Very recently, in \cite{Barker22} Barker shows how to quantify the main result in \cite{ESS}. He proves for $\delta$ small enough,
\begin{equation*}
    \limsup_{t \uparrow T^* } \frac {\| {\bf u} (\cdot ,t) \|_{L^3(B_{x_0} ( \delta) \, )}} { \left( \text{log log log}  \frac{1} {(T^* -t)^\frac{1}{4} } \right)^{\frac{1} {1129}} }  = + \infty,
\end{equation*}
improving thus \eqref{tao}.
All these advances have a ``supercritical" appearance however a less optimistic and more accurate reading would be : these advances are slightly supercritical \textbf{with respect to the $L^3$ norm}. With respect to the $L^{3,\infty}$ norm, these results have not a slightly supercritical relationship. Although $L^{3,\infty}$ is still a critical space, analogous results to those cited above, using $L^{3,\infty}$ instead of $L^3$, have not been demonstrated, except when conditions of smallness or axisymmetry are imposed. Thus, $L^{3,\infty}$ type conditions for the profiles, and $L^\infty((0,T),L^{3,\infty})$ type conditions in space-time, represent a big obstacle for the understanding of the Navier-Stokes equations.

To reinforce our opinion on the local well behavior of Navier-Stokes solutions, considering $L^3$ norms in the hypothesis, we will revisit some results.

\section{Removing smallness of the critical $L^3$ norm in localized smoothing}

In which concerns localized smoothing, we revisit Theorem 1 in \cite{TB_CP20} where a local smallness hypothesis of $L^3$ nature appear. To remove this smallness hypothesis the price to be paid is roughly the non-quantification of this result. Let us consider the definitions and notations in \cite{TB_CP20} (almost standard), in particular consider the constant $\gamma_{univ}$ described in Theorem 1 in \cite{TB_CP20}. Our main result reads as follow :
\begin{Theorem} 
\label{cor1}
Let ${\bf u} _0$ be a divergence free vector field belonging to $L^2_{\rm uloc}(\mathbb R ^3)$. 
Consider a compact set $K \subset  \mathbb{R}^3$.
We suppose that for each $x \in K$, there exists $\lambda_{x} > 0 $, such that $\int_{B_{2\lambda_{x}} (x)} |{\bf u}_0|^3 dy < \gamma_{\rm univ}^3$ and $ \|\lambda_{x} {\bf u}_0 (\lambda_{x} \cdot + x) \|_{L^2_{\rm uloc}} < + \infty$.

\noindent
Suppose $\bf u$ to be a local energy solution of (NS) with initial data ${\bf u}_0$ (Definition 16 in \cite{TB_CP20}). Then, there exists $T= T_{{\bf u}_0 , K}$ such that $ {\bf u} \in L^{\infty} (K \times (0, T)) .$
\end{Theorem}

\label{sec:intro}

{\bf Proof }
We will apply a scaling argument. 
Let $x \in K$. We take $\lambda_{x}$ such that 
\begin{equation*}
\label{hipT1}
    \int_{B_{2 \lambda_{x} }(x  ) } |{\bf u} _0 |^3 dy < \gamma_{univ}^3
\end{equation*}
and 
$ \| \lambda_{x} {\bf u}_0 (\lambda_{x} \cdot + x) \|_{L^2_{\rm uloc}} < + \infty. $
We have 
\begin{equation*}
    \int_{B_{2 \lambda_{x} }(x  ) } |{\bf u} _0 (y)|^3 dy = \int_{B_2(0  ) } |\lambda_{x} {\bf u} _0 (\lambda_{x} y + x) |^3 dy
\end{equation*}
and ${\bf u}_{\lambda_{x},x}(t,y)   :=   {\bf} \lambda_{x} {\bf u} (\lambda_{x}^2t , \lambda_{x} y + x )$ is a solution of the Navier-Stokes equations with initial data ${\bf u}_{0,\lambda_{x},x}(y) = \lambda_{x} {\bf u} _0 (\lambda_{x} y + x)$.

Let us denote $M_{{\bf u}_0, x}:= \| {\bf u}_{0,\lambda_{x},x} \|_{L^2_{\rm uloc}}$.
By Theorem 1 in \cite{TB_CP20}, we conclude that there exists $S^* (M_{{\bf u}_0, x})$ such that 
\begin{equation*}
    {\bf u}_{\lambda_{x},x} \in L^{\infty} (B_{\frac{1}{3}}(0) \times (\beta, S^* (M_{{\bf u}_0, x})))
\end{equation*}
for all $\beta \in (0, S^* (M_{{\bf u}_0, x}))$.
It follows that 
\begin{equation*}
    {\bf u} \in L^{\infty} (B_{\frac{\lambda_{x}}{3}}(x) \times (\beta, \lambda_{x}^2 S^* (M_{{\bf u}_0, x})))
\end{equation*}
for all $\beta \in (0, \lambda_{x}^2  S^* (M_{{\bf u}_0, x}))$.

As $\{  B_{\frac{\lambda_{x}}{3}}(x) ; x \in K \}$ cover the compact $K$, we can take a finite family $\{ x_{(1)}, x_{(2)},..., x_{(k)}\} \subset  \mathbb{R}^3$ such that 
\begin{equation*}
    \left\{ B_{\frac{\lambda_{x_{(i)}}}{3}}(x_{(i)}); i \in \{1,...,k\} \right\}
\end{equation*}
cover $K$. Then, we define 
$$T=T({{\bf u}_0,K})= \min_{i \in \{1,...,k\}}  \lambda_{x_{(i)}}^2 S^* (M_{{\bf u}_0, x_{(i)}})$$  
so that $ {\bf u} \in L^{\infty} (K \times (\beta,T)) $,
for all $\beta \in (0, T)$.

\begin{Remark}
If $K$ is not compact, to make this proof works, when we take for each $x \in K$ the positive number $\lambda_{x}$ such that
\begin{equation*}
    \int_{B_{2 \lambda_{x} }(x  ) } |{\bf u} _0 |^3 dy < \gamma_{univ}^3,
\end{equation*}
we need at the same time $\inf_{x \in K } \lambda_{x}^2 S^*(M_{{\bf u}_0, x})$ to be a positive number. So, we can impose as hypothesis $ \lambda_x $ to be bounded from below and $M_{{\bf u}_0, x}$ to be bounded from above. Then the following Corollary follows immediately from the proof of Theorem \ref{cor1}:
\end{Remark}

\begin{Corollary} 
\label{cor2}
Let ${\bf u} _0$ be a divergence free vector field belonging to $L^2_{\rm uloc}(\mathbb R ^3)$. We assume that
\begin{equation}
\label{hipMorrey}
    \sup_{x\in \mathbb R ^3} \sup_{\lambda >0 } \|\lambda{\bf u}_0 (\lambda \cdot + x) \|_{L^2_{\rm uloc}} < + \infty
\end{equation}
and there exists $\underline{\lambda} > 0$ such that for each $x \in \mathbb R ^3$, there exists $\lambda_{x} > \underline{\lambda} $, such that 
\begin{equation*}
    \int_{B_{2\lambda_{x}} (x)} |{\bf u}_0|^3 dy < \gamma_{\rm univ}^3.
\end{equation*}

\noindent
Suppose $\bf u$ to be a local energy solution of the Navier-Stokes equations with initial data ${\bf u}_0$. Then, there exists $T= T_{{\bf u}_0}$ such that $ {\bf u} \in L^{\infty} (\mathbb R ^3 \times (0, T)) $.
\end{Corollary}

\begin{Remark}
   We observe that \eqref{hipMorrey} is satisfied by the functions in the critical Morrey space $\dot{M}^{2,3}$, which contains $L^{3,\infty}$, where the norm is given by
   \begin{equation*}
       \|f\|_{\dot{M}^{2,3}}= \sup_{x \in \mathbb R^3} \sup_{R>0} \left( \frac{1}{R} \int_{|y-x|<R} |f|^2 dy\right)^{1/2}.
   \end{equation*}
This remark is also used for example in \cite{LRibe}. Observe that we could use Theorem 1.1 in \cite{KMT_21} in order to prove Theorem \ref{th2}.

\end{Remark}

\section{Removing backward discretely selfsimilar singularities}

In this section we want to show how scaling arguments permit to remove easily $L^3$ profiles in backward discretely selfsimilar  singularities. We will remove even $L^{3,q}(\mathbb R ^3)$ profiles with $1 \leq q <+ \infty $, where we consider the Lorentz spaces $L^{p,q}$, for $1\leq p,q < + \infty$, on measurable sets $\Omega\subset \mathbb R ^3$, defined by the quantity
\begin{equation*}
    \| f \|_{L^{p,q}( \Omega)} = \left( p \int_0^{+\infty} \alpha^q \, d_{f, \Omega}(\alpha) ^{q/p}  \, 
 \frac{d\alpha}{\alpha} \right)^{1/q}
\end{equation*}
where, 
\begin{equation*}
    d_{f, \Omega}(\alpha) = \mu ( \{ 
x \in \Omega \, : \, |f(x)| >\alpha \} ),
\end{equation*}
denoting $\mu$ the Lebesgue measure. That permits to precise some statements of recent results as Theorem 2 in \cite{CW17rem} or Corollary 1.1 in \cite{TB_CP21}. 

\noindent
We start by remember the definition of the $\lambda$-discretely self-similarity (see \cite{CW18,PF_PG}): 
\begin{Definition} 
    Let $\lambda > 1$.
\begin{itemize}
\item[$\bullet$] A vector field  $\vu_0 \in L^2_{\rm loc}(\mathbb{R}^3)$ is $\lambda$-discretely self-similar ($\lambda$-DSS)
if  $\lambda\vu_0(\lambda x)= \vu_0(x)$.
\item[$\bullet$] A time dependent vector field   $\vu \in L^2_{\rm loc}((-\infty,0)\times\mathbb{R}^3)$  is $\lambda$-DSS if  $\lambda \vu(\lambda^2 t,\lambda x)=\vu(t,x)$.
\item[$\bullet$] A pressure  $p \in L^{1}_{\rm loc}((-\infty,0)\times\mathbb{R}^3) $ is $\lambda$-discretely self-similar if  $ \lambda^2 p(\lambda^2 t, \lambda x) = p(t,x)$.
\end{itemize}
\end{Definition} 

\noindent
In the next statement, the notion of suitable Leray solution admits the standard definitions, including infinite energy solutions, we can consider for example the definition in \cite{TB_CP21}, \cite{CW17rem}, \cite{PF_PG}, \cite{LR02} or \cite{LR16}.
\begin{Theorem} 
\label{th2} 
Let $\lambda \in (1,+ \infty)$ and let ${\bf u}_{0}$ be a non trivial divergence-free vector field which is $\lambda$- DSS.  
Suppose, there exists a backward $\lambda$-DSS suitable Leray solution ${\bf u}$ of (NS) on $ (-\infty,0] \times \mathbb{R}^3 $.
Then, for all $t_0 \in (-\infty, 0) $, we have
\begin{equation}
\label{selfcrit}
    \liminf_{n \to \infty} \int_{B_1 \setminus B_{\lambda^{-1}}} |\lambda^n \vu(\lambda^n x,t_0) | >0.
\end{equation}
\end{Theorem}

\begin{Remark}
Observe that condition \eqref{selfcrit} implies the profile in negative times can't belong to $L^{3,q}$, with $1 \leq q <\infty$. In fact, by change of variables and by H\"older inequality (for Lorentz spaces),
\begin{align}
    \int_{B_1 \setminus B_{\lambda^{-1}} } |\lambda^n \vu(\lambda^n x,t_0) | dx &= \frac{1}{\lambda^{2n}} \int_{B_{\lambda^n} \setminus B_{\lambda^{n-1}} } | \vu( y,t_0) | dy \\
    & \leq \frac{1}{\lambda^{2n}}\| \mathds{1}_{B_{\lambda^n} \setminus B_{\lambda^{n-1}}} \|_{L^{\frac{3}{2}, \frac{q}{q-1}}}  \| \vu( \cdot ,t_0) \|_{L^{3,q}(B_{\lambda^n} \setminus B_{\lambda^{n-1}} )},
\end{align}
and since $\| \mathds{1}_{B_{\lambda^n} \setminus B_{\lambda^{n-1}}} \|_{L^{\frac{3}{2}, \frac{q}{q-1}}} \leq C \lambda^{2n}$, we find
\begin{align}
    \liminf_n \int_{B_1 \setminus B_{\lambda^{-1}} } |\lambda^n \vu(\lambda^n x,t_0) | \leq C \liminf_n  \| \vu( \cdot ,t_0) \|_{L^{3,q}(B_{\lambda^n} \setminus B_{\lambda^{n-1}} )}.
\end{align}
By other hand, since $\vu(t_0) \in L^{3,q}$, with $1 \leq q<\infty$, using the definition of this Lorentz space and the monotone convergence, we can verify that $ \lim_{n \to \infty} \| \vu( \cdot ,t_0) \|_{L^{3,q}(\mathbb R^3 \setminus B_{\lambda^{n-1}} )}=0$.

In terms of $L^2$ norms, we have \eqref{selfcrit} also implies 
\begin{equation*}
    \liminf_{R \to \infty} \frac{1}{R} \int_{B_R} |\vu(\cdot, t_0)|^2 >0.
\end{equation*}
In fact, by the H\"older inequality
\begin{align*}
    \liminf_{n \to \infty} \int_{B_1 \setminus B_{\lambda^{-1}} } |\lambda^n \vu(\lambda^n x,t_0) | = \liminf_{n \to \infty} \frac{1}{\lambda^{2n}} \int_{B_{\lambda^{n}} \setminus B_{\lambda^{n-1}}} | \vu( x,t_0) | \leq \liminf_{n \to \infty} \frac{c}{\lambda^{n}} \int_{B_{\lambda^{n}} \setminus B_{\lambda^{n-1}}} | \vu( x,t_0) |^2.
\end{align*}

\end{Remark}

{\bf Proof of Theorem \ref{th2} :}
Suppose there exists $t_0 \in (-\infty, 0) $ such that
\begin{equation}
\label{selfcrits}
    \liminf_{n \to \infty} \int_{B_1 \setminus B_{\lambda^{-1}}} |\lambda^n \vu(\lambda^n x,t_0) | =0.
\end{equation}
We will demonstrate $\vu(\cdot,0)=0$. It is sufficient to show that for almost all $x \in B_1 \setminus B_{\lambda^{-1}}$, $\vu(x,0)= \lim \vu(x,t_n)=0$. Then, by discretely selfsimilarity of the initial profile we obtain $\vu(\cdot,0)=0$.

\noindent
By interior estimates of regularity (see \cite{Serbook14} and \cite{Tsaibook17}) we know 
\begin{equation}
\label{cl1}
    \vu \in \mathcal{C}([-1,0],L^1(B_1)) .
\end{equation}
Consider the sequence $t_n = t_0 / \lambda^{2n}$. By $\lambda$-DSS property we get 
\begin{equation}
\label{seq}
    \lambda^n \vu(\lambda^n x,t_0)= \vu(x, \frac{1}{\lambda^{2n}} t_0).
\end{equation}
By \eqref{seq} and \eqref{selfcrits} we can extract a subsequence of $\vu(\cdot, \frac{1}{\lambda^{2n}} t_0)$ converging to zero pointwise. 
We so get, up to a subsequence, for almost all $x \in B_1 \setminus B_{\lambda^{-1}}$, $\vu(x,0)= \lim \vu(x,t_n)=0$ by \eqref{cl1}.

\begin{Corollary}
Suppose that $\vu$ defined on $(-\infty,0] \times \mathbb R ^3$ is a suitable backward discretely selfsimilar solution such that 
\begin{equation*}
    \| \vu \|_{L^\infty((-1,0),L^{3,\infty})} < +\infty
\end{equation*}
and there exists $t_0 <0$ such that $\vu(t_0) \in L^{3,q}$, with $3 \leq q < +\infty$, then
$\vu=0$.
\end{Corollary}
\noindent
It is a direct consequence of Theorem \ref{th2} and Theorem 4.1 in \cite{AB19}
(a Liouville type theorem).

\small{DATA AVAILABILITY STATEMENT} : \small{Data sharing isn't applicable to this article as no datasets were generated or analyzed during
the current study.}
\section*{\small{ACKNOWLEDGMENT}}
\noindent
\small{Pedro Fernández is partially supported by the Project PEPS JCJC 2023–UMR 8088 (AGM), Insmi.}


\begin{thebibliography}{HD}

\bibitem{AB19} D. Albritton and T. Barker,
    \textit{On Local Type I Singularities of the Navier–Stokes Equations and Liouville Theorems}, J. Math. Fluid Mech. (2019).

\bibitem{BFP} T. Barker, P.-G Fern\'andez-Dalgo and Christophe Prange.
    \textit{Blow-up of dynamically restricted critical norms near a potential Navier-Stokes singularity}, Math. Ann. (2023).

\bibitem{Barker22} T. Barker,
    \textit{Localized quantitative estimates and potential blow-up rates for the Navier-Stokes equations}, To appear in Siam Journal on Mathematical Analysis.

\bibitem{TB_CP20} T. Barker and C. Prange,
    \textit{Localized Smoothing for the Navier–Stokes Equations and Concentration of
    Critical Norms Near Singularities}, Arch. Ration. Mech. Anal., 236(3):1487–1541, (2020).

\bibitem{TB_CP21} T. Barker and C. Prange,
    \textit{Quantitative regularity for the Navier-Stokes
    equations via spatial concentration.}, Comm. Math. Phys., 385(2):717–792, (2021).

\bibitem{CW17rem} 	D. Chae and J. Wolf,  
    \emph{Removing discretely self-similar singularities for the 3D Navier–Stokes equations}, Communications in Partial Differential Equations, (2017).

\bibitem{CW18} 	D. Chae and J. Wolf,  
    \emph{Existence of discretely self-similar solutions to the Navier-Stokes
    equations for initial value   in $L^2_{\rm loc}(\mathbb{R}^3)$}, Ann. Inst. H. Poincar\'e  Anal. Non Lin\'eaire  {35}, 1019--1039 (2018).	

\bibitem{ESS} L. Escauriaza, G. Seregin, V. Šverák, 
    \emph{$L_{3,\infty}$-Solutions to the Navier-Stokes Equations and Backward
    Uniqueness}, Russ. Math. Sur. 58(2), 211–250 (2003)

\bibitem{PF_PG} P.-G Fern\'andez-Dalgo and P.-G.  Lemari\'e--Rieusset, 
    \emph{Weak Solutions for Navier–Stokes Equations with Initial Data in Weighted L2 Spaces.}, Archive for Rational Mechanics $\&$ Analysis (2020)

\bibitem{FK_64} H. Fujita and T. Kato,
\textit{On the Navier-Stokes initial value problem. I}, Archive for Rational Mechanics and Analysis volume 16, pages 269–315 (1964)

\bibitem{KMT_21} K. Kang, H. Miura and T. Tsai,
    \textit{An  $\epsilon$-regularity criterion and estimates of the regular set for Navier-Stokes flows in terms of initial data}, Pure Appl. Analysis 3 567-594 (2021).

\bibitem{La67} 
O. Ladyzhenskaya, \emph{On uniqueness and smoothness of generalized solutions to the Navier-Stokes equations}, Zapiski Nauchn. Seminar. POMI,
5, 169–185 (1967).

\bibitem{LRibe} P.-G. Lemarié-Rieusset,
\textit{The Navier-Stokes equations in the critical Morrey-Campanato space},
Rev. Mat. Iberoamericana 23(3): 897-930, 2007
	
\bibitem{LR02}  P.-G. Lemari\'e-Rieusset, 
    \emph{Recent developments in the Navier-Stokes problem}, CRC Press, (2002).

\bibitem{LR16} P.-G.  Lemari\'e-Rieusset,       
    \emph{The Navier-Stokes problem in the 21st century}, Chapman \& Hall/CRC, (2016).

\bibitem{Le34} J. Leray, 
    \emph{Essai sur le mouvement d'un fluide visqueux emplissant l'espace}, Acta Math. {63} 193--248 (1934).

\bibitem{P59}
    G. Prodi, \emph{Un teorema di unicit\'a per le equazioni de Navier-Stokes}, Ann. Mat. Pura Appl., 48:173-182, (1959).

\bibitem{Ser12}
    G. Seregin, \emph{A certain necessary condition of potential blow up for Navier-Stokes equations}, Communications in Mathematical Physics, 312(3):833–845, (2012).

\bibitem{Serbook14}
    G. Seregin, \emph{Lecture notes on regularity theory for the Navier-Stokes equations}, World Scientific, (2014).

\bibitem{S62} J. Serrin, 
    \emph{On the interior regularity of weak solutions of the Navier-Stokes equations}, Arch. Ration. Mech. Anal., 9 pp. 187–195 (1962).

\bibitem{Tao19} T. Tao, 
    \emph{Quantitative bounds for critically bounded solutions to the Navier-Stokes equations}, In A. Kechris, N. Makarov, D. Ramakrishnan, and X. Zhu, editors, Nine Mathematical Challenges: An Elucidation, volume 104. American Mathematical Society, (2021).

\bibitem{Tsaibook17} T.-P. Tsai,
    \textit{Lectures on Navier-Stokes equations}, American Mathematical Society, Graduate Studies in Mathematics 192, (2017).

\end{thebibliography}
\end{document}